%%%%%%%%%%%%%%%%%%%%%%%%%%%%%%%%%%%%%%%%%%%%%%%%%%%%%%%%%
%	 I. Frenkel, A. Kirillov Jr, and A. Varchenko 			%
%        Canonical basis and homology of local systems  %
%                                                       %
%%%%%%%%%%%%%%%%%%%%%%%%%%%%%%%%%%%%%%%%%%%%%%%%%%%%%%%%%
% This file should be processed by AmsTeX version 2.1%%%%
%%%%%%%%%%%%%%%%%%%%%%%%%%%%%%%%%%%%%%%%%%%%%%%%%%%%%%%%%

\input epsf
\input amstex
\loadbold
\magnification=1200
\documentstyle{amsppt}
\NoBlackBoxes
\def\Z{\Bbb Z}
\def\C{\Bbb C}
\def\R{\Bbb R}
\def\Q{\Bbb Q}
\def\a{\alpha}
\def\l{\lambda}
\def\bl{{\boldsymbol \lambda}}
\def\zt{\tilde z}
\def\g{\frak g}

\def\Id{\operatorname{Id}}
\def\Hom{\operatorname{Hom}}

\def\sltwo{\frak s\frak l _2 }

\def\U{U_q \frak s\frak l _2 }

\def\<{\langle}
\def\>{\rangle}

\def\i{\text{\rm i}}
\def\o{\otimes}
\def\half{{\frac{1}{2}}}

\topmatter
\title Canonical basis and homology of local systems 
\endtitle
\date March 4, 1997 \enddate
\author Igor Frenkel, Alexander Kirillov, Jr. \\
 and  Alexander Varchenko\endauthor
\leftheadtext{Frenkel, Kirillov and Varchenko}
\endtopmatter

\document

\head Introduction\endhead

Let $\U$ be the quantum group corresponding to the Lie algebra
$\sltwo$, and let $M_\l$ be the Verma module over the $\U$ with
the highest weight $\l$. It follows from the results of 
\cite{V1} that one can identify the weight space
$M[\mu]=(M_{\l_1}\o \dots \o M_{\l_n})[\mu], \mu=\sum \l_i- 2l$ with a
suitable homology space of the configuration space $X_{n,l}=\{\bold
x\in \C^l|x_i\ne x_j, x_i\ne 0, x_i\ne z_k, k=1,
\dots, n\}$ with coefficients in a certain one-dimensional local
system.  This result plays a crucial role in the construction of quantum
group symmetries in conformal field theory, since this homology space
naturally appears in the integral formulas for the solutions of
Knizhnik-Zamolodchikov equations (see \cite{SV1, SV2}).

On the other hand, for any simple Lie algebra $\g$, Lusztig has
defined a remarkable basis in a tensor product of irreducible
finite-dimensional representations of $U_q\g$, which he called ``the
canonical basis''(see \cite{L, Chapter 27}). This basis is defined by
two conditions:

(1) $\Psi(b)=b$, where $\Psi$ is a certain antilinear involution
involving the $R$-matrix.

(2) The canonical basis is related with the usual basis, given by 
the tensor product of Poincare-Birkhoff-Witt bases,  by a matrix, all
off-diagonal entries of which belong to $q^{-1}\Z[q^{-1}]$.

This basis generalizes the canonical basis in one irreducible
finite-dimensional module over $U_q\g$, also introduced by Lusztig, and
has a number of remarkable properties (see \cite{L}).

Note that the canonical basis in a tensor product is non-trivial even
for $\U$, since the definition of $\Psi$ involves the quantum
$R$-matrix.  For this reason, the first condition is usually much more
difficult to verify than the second one.

The goal of this paper is to combine these two results and give a
geometric construction of the canonical basis in terms of the
homologies of local systems (for technical reasons, we are
constructing the basis dual to the canonical basis rather than the
canonical basis itself). In this paper, we only consider $\g=\sltwo$
and assume that $q$ is not a root of unity. In this case we are able
to give an explicit construction of the dual canonical basis in every
weight subspace of a tensor product of irreducible finite-dimensional
$\U$-modules (see Theorem~4.3).  This answer is especially simple if
all the highest weights are large enough compared to the level $l=(\sum
\l_i -\mu)/2$ of the considered weight subspace. In this case, the
basis dual to the canonical basis is given (up to some simple factor)
by the bounded connected components in the complement to certain
hyperplanes in $\R^l$ (Corollary~4.7). Our construction is quite
parallel to the algebraic constructions in \cite{FK}.

Notice that in \cite{V2} a geometric construction of the crystal base in the 
space of singular vectors of $M[\mu]$ is given in terms
of the same local system and critical points of the associated
multivalued holomorphic function. It would be interesting to establish
a direct connection between these two constructions.

Generalization of the results of this paper to arbitrary simple Lie
algebras is more complicated; for example, even for $\frak{sl}_3$ the
natural generalization of the approach in \cite{FK} fails to produce
the canonical basis (see \cite{KK}). We plan to address these
questions in forthcoming papers.

\head 1. $\U$: notations \endhead

Let the quantum group $\U$ be the Hopf algebra over the
field $\C(q^{\pm\half})$ with generators $e,f, q^{\pm h}$ and
commutation relations 
$$\gathered
q^h e =q^2 eq^h,\\
q^h f =q^{-2} fq^h,\\
[e,f]=\frac{q^h-q^{-h}}{q-q^{-1}}.\endgathered
\tag 1.1$$

The comultiplication and antipode are given by 
$$\gathered 
\Delta e= e \o q^{h/2} + q^{-h/2}\o e,\\
\Delta f= f \o q^{h/2} + q^{-h/2}\o f,\\
\Delta q^h=q^h\o q^h,\\
Se=-qe,\qquad Sf=-q^{-1}f,\qquad S q^h=q^{-h}.
\endgathered\tag 1.2$$

This coincides with the definition in \cite{V1, 4.1} if we make the
following substitutions: $q\mapsto q^\half, e\mapsto e/\sqrt{q^\half
-q^{-\half}}, f\mapsto f/\sqrt{q^\half-q^{-\half}}$. 

We will also use another generators of the same algebra, namely:

$$E=q^{h/2}e,\qquad F=fq^{-h/2}. \tag 1.3$$
They satisfy the same commutation relations (1.1) as $e,f$, and the
comultiplication is given by 
$$\gathered
\Delta E=E\o q^h+ 1\o E,\\
\Delta F=F\o 1+ q^{-h}\o F.
\endgathered 
\tag 1.4
$$

Also, it is useful to note that 

$$\gathered
E^m=q^{m(m+1)/2 } e^m q^{mh/2},\\
F^m=q^{m(m-1)/2 } f^m q^{-mh/2}.
\endgathered
\tag 1.5
$$

The universal R-matrix for $\U$ is given by

$$\gathered
\Cal R=C\Theta,\quad C=q^{h\o h/2},\\
{\aligned
	\Theta&=\sum_{k\ge 0} q^{-k(k+1)/2}
		   \frac{(q-q^{-1})^k}{[k]!}q^{kh/2}e^k\o q^{-kh/2}f^k\\
	 &=\sum_{k\ge 0} q^{k(k-1)/2}
		\frac{(q-q^{-1})^k}{[k]!} E^k\o F^k.
\endaligned}
\endgathered
\tag 1.6
$$

Here, as usual, $[n]=\frac{q^n-q^{-n}}{q-q^{-1}}$, and
$[n]!=[1][2]\dots [n]$. 

For every pair $V,W$ of finite-dimensional representations, the  $R$-matrix
gives rise to the commutativity isomorphism $\check \Cal R=P \Cal R:V\o W\to
W\o V$, where $P$ is the permutation: $P(v\o w)=w\o v $.

We shall only consider $\U$-modules with  weight
decomposition: $V=\bigoplus_{\l\in \C} V[\l]$, and $q^h|_{V[\l]}=q^\l\Id$
(for  a non-integer $\l$, this requires an appropriate extension of the
field of scalars).

We denote by $M_\l$ the Verma module with highest weight $\l$. For $\l\in
\Z_+$ this module has a finite-dimensional quotient: 

$$V_\l=M_\l/(f^{\l+1} v_\l).$$

Then $V_\l$ is an irreducible finite-dimensional module of dimension
$\l+1$, and every  irreducible finite-dimensional  module with weight
decomposition is isomorphic to one of $V_\l$. 

Finally, define an algebra anti-automorphism  $\tau:\U\to \U$ 
by 

$$\tau(e)=f, \quad \tau(f)=e, \quad \tau (q^h)=q^h,
	\quad \tau(ab)=\tau(b)\tau (a).\tag 1.7$$
Then $\tau$ is a coalgebra automorphism:
$(\tau\o\tau)\Delta(x)=\Delta(\tau(x))$,   and $\tau (\Cal R)=\Cal
R^{21}$. From now on, we will also denote by $\tau$ the map
$\tau\o\dots \o \tau: (\U)^{\o n}\to (\U)^{\o n}$. 

For every module $M$ let the contragredient module $M^c$ be the
restricted dual to $M$ with the action of $\U$ given by 

$$\<gv^*, v\> = \<v^*, \tau(g)v\>, \qquad v\in M, 
		\quad v^*\in M^c,\quad g\in \U.
\tag 1.8$$
Note that we have a canonical isomorphism $(M_1\o M_2)^c\simeq M_1^c\o
M_2^c$ and that for $\l\in \Z_+$, $V_\l^c\simeq V_\l$. Also, for every
$\l$ we have a canonical morphism $M_\l\to M_\l^c$, and the image of
this morphism is exactly the irreducible highest-weight module $L_\l$;
in particular, for $\l\in \Z_+$, the image is $V_\l$. Thus, we have an
embedding 

$$V_\l\subset M_\l^c. \tag 1.9$$

Combining this embedding with the canonical pairing $M_\l\o M_\l^c\to
\C(q^{\pm\half})$, we get a non-singular bilinear form on $V_\l$, which is usually
called the contragredient (or Shapovalov) form.

All the theory above can be as well developed when $q$ is a complex
number rather than a formal variable, provided that $q$ is not a root
of unity. We will use the same notations in this case.

\remark{Remark 1.1} Our notations differ slightly from those of
Lusztig. They are related as follows:

$$
\alignat2
&\text{Ours} \quad && \text{Lusztig's}\\
& q && v\\
& E && F\\
&F && E\\
&q^{h} && K^{-1}\\
&\Theta && \bar\Theta
\endalignat
$$
(in the last line, $\overline{\phantom{T}}$ is the bar involution,
see (3.1)). 
\endremark

\head 2. Homology of local systems\endhead

Fix the following data:

-- a positive integer $n$

-- a collection of weights $\l_1, \dots, \l_n\in\C$

-- points $z_1,\dots, z_n\in \R, 0<z_1<\dots<z_n$ 

-- an integer $l\in \Z_+$

-- a  number $\kappa\in \R\setminus \Q$

We denote $\l=\l_1+\dots +\l_n$ and 

$$q=e^{\pi\i/\kappa}.
 $$
Note that the number $q$ defined by the formula above is not a root of
unity.

Let $\Cal C\subset \C^l=\{(x_1, \dots, x_l)\}$ be the configuration of
hyperplanes defined by the equations

$$\gathered 
x_i=x_j, \quad i,j=1\dots l,\\
x_i=z_k, \quad i=1,\dots, l, \quad k=1,\dots, n,\\
x_i=0, \quad i=1,\dots, l.
\endgathered
\tag 2.1
$$

Let $Y=\C^l\setminus \Cal C$, and let $\Cal S$ be the one-dimensional
complex local system over $Y$ such that its flat sections are
$s(x)=const(\text{univalent branch of }\psi(x))$, where
$const \in \C$ and 

$$\psi=\prod_{i<j} (x_i-x_j)^{2/\kappa}
	\prod_{i,p} (x_i-z_p)^{-\l_p/\kappa}
	\prod_{p<q} (z_p-z_q)^{\l_p\l_q/2\kappa}.
\tag 2.2
$$

Let $\Sigma_l$ be the symmetric group acting on $\C^l$ by permutation
of coordinates; this action preserves $\Cal C$ and $\Cal S$. We will be
interested in the relative homology group $H_l(\C^l; \Cal C; \Cal S)$
with coefficients in $\Cal S$. We have a natural action of
$\Sigma_l$ on this homology space. Let us for brevity denote by $H_l$
the antisymmetric part of the homology: 

$$H_l=H_l^{-\Sigma_l}(\C^l; \Cal C; \Cal S).\tag 2.3$$

Sometimes we will also write $H_l(0, z_1,\dots, z_n; 0,
\l_1,\dots,\l_n)$ to emphasize the dependence on $z_i, \l_i$ (the
first zero is to remind of the point $0$ also used in the definition;
it will be needed later when we define slightly more general type of
homologies).  In fact, the space $H_l(0,\bold z;0, \bl )$ is defined
for any $\bold z\in X_n=\{(z_1, \dots, z_n)\in \C^n|z_i\ne z_j, z_i\ne
0\}$, and $H_l(0,\bold z;0, \bl)$ form a vector bundle over
$X_n$. Moreover, this bundle has a natural flat connection
(Gauss-Manin connection, see \cite{V1}). However, unless otherwise
noted, we will only use $\bold z$ such that $0<z_1<\dots <z_n$.

We will frequently use pictures to represent cycles in $H_l(\C^l; \Cal
C; \Cal S)$. For example, we will represent the cycle $z_1\le x_1\le
x_3\le x_2\le z_2$ by the picture

$$\epsfbox{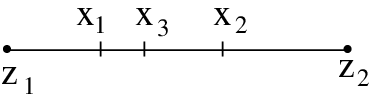}
$$

We will choose a section $s$ of the local system over such a cycle as
follows. For each pair of indices $i,j$ define 

$$Br((x_i-x_j)^\a)=\cases e^{\a\text{Log}(x_i - x_j)}, & Re (x_i)>Re (x_j),\\
			e^{\a\text{Log}(x_j- x_i)}, & Re (x_i)<Re (x_j),
		\endcases
$$
where $\text{Log}(x)$ is the main branch of the logarithm defined for
$Re\ x> 0$ by the condition $\text{Log}\ x\in \R_+$ for $x>0$. This
definition is chosen so that for $x_i, x_j\in \R_+, \a \in \R$, $Br
(x_i-x_j)^\a\in \R_+$.  Similarly, define

$$Br (\psi)= \prod Br(x_i-x_j)^{2/\kappa}
	\prod Br(x_i-z_k)^{-\l_k/\kappa}
	\prod Br(z_i-z_j)^{\l_i\l_j/2\kappa}.
$$

 One easily sees that if a region $D\subset\C^l$ is such that
$x_i-x_j\notin \i\R, x_i-z_k\notin \i\R$ on $D$, then $Br(\psi)$ is a
section of the local system $\Cal S$ on $D$. In particular, if we have a
simplex $c\subset \C^l$ such that the order of $Re\ x_i, Re\ z_i$ is
fixed on $c$ then $Br (\psi)$ defines a section of $\Cal S$ over $c$.  In
\cite{V1}, this section is called ``the positive branch of $\psi$''.

  We will also use more elaborate $l$-dimensional cycles, such as the
one used in formula (2.6) below. It will be convenient to give the
following definition; it may seem complicated but is in fact quite
natural.

\definition{Definition 2.1}
Assume that we have fixed the data $n, \bold z, \bl, l, \kappa$ as at
the beginning of this section. Let us additionally assume that $Re\
z_i\ne Re\ z_j, Re\ z_i\ne 0$.  A ``combinatorial cycle'' $C$ is a
figure in $\C$ consisting of a finite number of intervals of smooth
curves \rom{(}``arcs''\rom{)} such that:

- these arcs do not intersect

- each arc can intersect a vertical line at most at one point 

- the  endpoints of each arc are either some of the specified points $z_i$ or $0$, and
the endpoints of the same arc can not coincide

- on each arc, there are some marked points so that the total
number of the marked points on all arcs is equal to $l$

- these marked points are labeled by $x_1, \dots, x_l$

\enddefinition

For every combinatorial cycle $C$ as above, we define the
corresponding relative cycle in $H_l(\C^l, \Cal C;\Z)$ by the
following rule. Consider the space $\R^l$ with the coordinates
$x_1,\dots, x_l$ and the standard orientation. Let $\Delta_C$ be the
subset in $\R^l$ given by the following conditions:

- $0\le x_i\le 1$

- if the points $x_i$ and $x_j$ are marked on the same arc of $C$, and
$x_i$ is to the left of $x_j$ then $x_i\le x_j$ for all points in
$\Delta_C$.

Obviously, $\Delta_C$ is a direct product of simplices. 

Now, let us number all the arcs of $C$ in an arbitrary way and let
$\gamma_k:[0,1]\to \C$ be a parameterization of the $k$-th arc such
that $\text{\rm Re}\ \gamma(t)$ is a strictly increasing
function. Define the map $\gamma:\Delta_C\to \C^l$ by

$$\gamma(x_1,\dots, x_l)=(\gamma_{k_1}(x_1),\dots,
\gamma_{k_l}(x_l)),$$
where $k_i$ is the index of the arc on which the point $x_i$ is
marked. This defines an $l$-dimensional chain in $\C^l$, which we will
also denote by $C$. Informally, this chain can be described by saying
that the points $x_i$ run along the corresponding arcs preserving the
order.

 Define a section $s$ of the local system $\Cal S$ over $C$ by
the rule

$$s(\bold x)=q^{\sum a_i} Br (\psi(\bold x)),$$
where $\bold x$ is such a point that all $x_i$ are close to the right
end of the corresponding arc, and the order of $Re\  x_i$ is chosen as
in the following picture:

$$\epsfbox{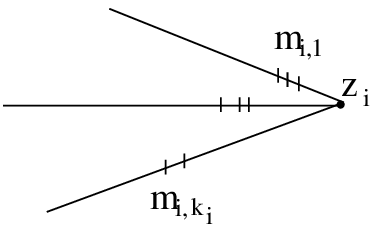}$$

and $a_i=-\sum_{1\le a<b\le k_i}m_{i,a}m_{i,b}b$

\proclaim{Lemma 2.2} Let $C$ be a combinatorial cycle. Then the pair
$(C, s)$ defines an element of the relative homology space
$H_l(\C^l, \Cal C; \Cal S)$, which only depends on the homotopy type
of the arcs and the order of the points placed on these arcs. These
cycles are flat sections of the  bundle $H_l(\C^l, \Cal C; \Cal S)$
over $X_n$ with respect to the Gauss-Manin connection.
\endproclaim

The reason for calling these cycles ``combinatorial'' is that they are
defined by a finite collection of data. Later we will show that such
cycles span the whole homology space $H_l$ (see Theorem~2.3). 

Note that the action of a permutation $\sigma\in \Sigma_l$ on the
combinatorial cycles is given by permuting the indices of the points
$x_i$ and multiplying by $(-1)^{|\sigma|}$ (the sign comes from the change
of orientation). Since we are only interested in the antisymmetric
part of the homology, from now on  we won't put
the labels $x_i$ on the points.  Instead, we will just indicate the
number of points on each arc, assuming  summation over all $l!$
possible labelings, which automatically gives an anti-symmetric cycle.

Let $M_\l$ be a Verma module over $\U$ with $q=e^{\pi\i/\kappa}$ and let
$v_\l$ be a highest weight vector. Vectors $F^{(m)}v_\l=\frac{F^m
v_\l}{[m]!}, m=0,1,\dots$ form a basis in $M_\l$. Denote by
$(F^{(m)}v_\l)^*$ the dual basis in $M_\l^c$. More generally, for any
$n$-tuple $\bold m=(m_1, \dots, m_n)\in \Z_+^n$ denote

$$F^{(\bold m)}=F^{(m_1)}v_{\l_1}\o \dots\o F^{(m_n)}v_{\l_n}\in
M_{\l_1}\o\dots \o M_{\l_n}.
\tag 2.4
$$
These monomials form a basis in $M_{\l_1}\o\dots\o M_{\l_n}$. Let 

$$(F^{(\bold m)})^*=(F^{(m_1)}v_{\l_1})^*\o \dots\o (F^{(m_n)}v_{\l_n})^*
\tag 2.5$$
be the dual basis in $M^c_{\l_1}\o\dots \o M^c_{\l_n}$. It is easy to
see that $(F^{(\bold m)})^*\in M^c_{\l_1}\o\dots \o
M^c_{\l_n}[\sum \l_i-2\sum m_i]$ (recall that for any $\U$-module $V$
we denote by $V[\mu]$ the subspace of vectors of weight $\mu$ in $V$).

\proclaim{Theorem 2.3} Assume that $z_i\in \R,
0<z_1<z_2\dots<z_n$. Define the map 
$\varphi_{\bold z}:M_{\l_1}^c\o \dots \o M^c_{\l_n}[\l-2l]\to H_l$ 
by 

$$(F^{(\bold m)})^*
	\mapsto [m_1]!\dots [m_n]! 
\vcenter{\epsfbox{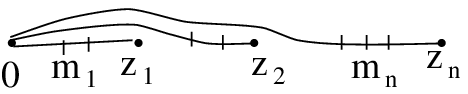}}
\tag 2.6
$$ 
where $l=\sum m_i$.

Then:
\roster 
\item Map \rom{(2.6)} is an isomorphism.

\item For any $i=1,\dots, n-1$, we have the following commutative diagram: 

$$\CD
M^c_{\l_1}\o\dots\o M^c_{\l_n}[\l-2l] @>{\varphi_{\bold z}}>> 
				H_l(0,\bold z; 0,\bl)\\
	@VV{\check \Cal R_i}V			@VV{T_i}V\\
M^c_{\l_1}\o\dots\o M^c_{\l_{i+1}}\o M^c_{\l_i}\o \dots\o M^c_{\l_n}[\l-2l]
	@>{\varphi_{s_i(\bold z)}}>>	H_l(0,s_i(\bold z), 0,s_i(\bl))
\endCD
$$
where  $s_i(\bold z)=(z_1, \dots, z_{i+1}, z_i, \dots, z_n), 
	s_i(\bl)=(\l_1, \dots, \l_{i+1}, \l_i, \dots, \l_n)$ 
and $T_i$ is the monodromy along the path shown on Figure~1. 

\endroster
\endproclaim

\midinsert
\centerline{\epsfbox{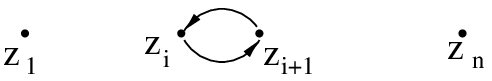}} 
\botcaption{Figure 1}
\endcaption
\endinsert

\demo{Proof} This theorem follows from the results of  \cite{V1}. For
the sake of completeness, we quote here the statement in the most
convenient for us form.  Let $z_0<\dots<z_n\in \R$, and let $\Cal
C(z_0, \dots, z_n)$ be a configuration of hyperplanes in $\C^l$, given
by the equations $x_i=x_j$, $x_i=z_k, k=0\dots n$.  Let us also fix
real numbers $\l_0,
\dots, \l_n$ (``weights'') and let $\Cal S$ be the one-dimensional
local system on $\C^l\setminus \Cal C(z_0,\dots, z_n)$ defined as in
(2.2), but with the indices $p,q$ running from $0$ to $n$ rather than
from $1$ to $n$.  Denote by $H_l(z_0,\dots, z_n; \l_0,\dots,
\l_n)=H_l^{-\Sigma_l}(\C^l,\Cal C(z_0, \dots, z_n),\Cal S)$ the
antisymmetric part of the relative homology space. Obviously, if we
let $z_0=0, \l_0=0$ then this space coincides with previously
defined $H_l= H_l(0, z_1,\dots, z_n;0,\l_1,\dots, \l_n)$.

Recall that
a vector $v$ in a $\U$-module $V$ is called singular if $e v=0$. We
denote by $V^{sing}$ the subspace of singular vectors in \cite{V1}. 
Then the following result is a corollary of Theorem~5.11.13 in
\cite{V1} (see also Figure~5.17).

\proclaim{Theorem 2.4 \rm(\cite{V1})} Let $z_0,\dots, z_n\in
\R,z_0<\dots<z_n$.   
There exists an isomorphism

$$\varphi:  \bigl(M_{\l_0}^c\o \dots \o
M_{\l_n}^c\bigr)^{sing}[\sum_0^n\l_i-2l]\simeq 
 H_l(z_0,\dots, z_n;\l_0,\dots, \l_n)
\tag 2.7$$

such that 

\roster\item 
 Let $C_{m_1, \dots, m_n}\in H_l(z_0,\dots, z_n;\l_0,\dots, \l_n)$ be
the cycle in the right-hand side of \rom{(2.6)} \rom{(}with $0$
replaced by $z_0$\rom{)}. Then

$$\varphi^{-1}(C_{m_1, \dots, m_n})=v_{\l_0}^*\o (F^{(m_1)}v_{\l_1})^*\o \dots 
\o (F^{(m_n)}v_{\l_n})^*+\dots,$$
where dots stand for a combination of monomials of the form 
$ (F^{(m_0)}v_{\l_0})^*\o \dots\o  (F^{(m_n)}v_{\l_n})^*$ with $m_0>0$.

\item For any $i=0,\dots, n-1$, 
we have a commutative diagram similar to that in Theorem~\rom{2.3}
part \rom{(2)}. 
\endroster
\endproclaim

\remark{Remark} 
In \cite{V1}, the basis $f^m=f^{m_1}v_{\l_1}\o\dots \o f^{m_n}v_{\l_n}$
is used instead of $F^{(m)}$. In these notations, formula (2.6)
 takes the form

$$(f^{m_1} v_{\l_1})^*\o \dots\o (f^{m_n} v_{\l_n})^*\mapsto
q^{a(\bold m)}\vcenter{\epsfbox{sltwo3.eps}}
$$ 
where 
$$a(\bold m)=\sum_{i=1}^n \frac{m_i(m_i-1-\l_i)}{2}.
$$
(compare with \cite{V1, 5.2.22})
\endremark

In order to get our Theorem~2.3 from Theorem ~2.4, we also need the
following lemma. 

\proclaim{Lemma 2.5} Let $V$ be an arbitrary highest weight module over
$\U$, and let $M_{\l_0}^c$ be the contragredient Verma module with the
highest weight $\l_0$. Then: for every homogeneous vector $v\in
V[\mu]$ there exists a unique vector $\tilde v\in (M_{\l_0}^c\o
V)^{sing}[\l_0+\mu]$ such that 

$$\tilde v=v^*_{\l_0}\o v+\dots,$$
where dots stand for a combination of vectors of the form
$(F^{(m)}v_{\l_0})^*\o v'$ with $m>0$. 

This gives us an isomorphism 

$$\aligned
V[\mu]&\simeq (M_{\l_0}^c\o V)^{sing}[\l_0+\mu]\\
v&\mapsto \tilde v
\endaligned
\tag 2.8
$$

\endproclaim

\demo{Proof of the Lemma}Follows from  the following identities:

$$\split
(M_{\l_0}^c\o V)^{sing}[\l_0+\mu]
   =&\Hom_{U_q\frak b}(\C_{\l_0+\mu},M_{\l_0}^c\o V) 
	=\Hom_{U_q\frak b}( (M_{\l_0}^c)^*\o \C_{\l_0+\mu}, V)\\
	=&\Hom_{U_q\frak b}((M^c_{-\mu})^*, V)
	  =V[\mu],
\endsplit
$$
where $U_q\frak b$ is the subalgebra in $\U$ generated by $q^h$ and
$e$ and $\C_{\l_0+\mu}$ is the one-dimensional module over $U_q\frak
b$ with the action given by $e=0, q^h=q^{\l_0+\mu}$. The last identity
follows from the fact that $(M^c_{-\mu})^*$ is nothing but the lowest
weight Verma module with lowest weight $\mu$ and thus is free over
$U_q\frak n^+$.\qed
\enddemo

Now we can easily prove our Theorem~2.3. Indeed, letting in
Theorem~2.4 $\l_0=0, z_0=0$, we get an isomorphism $H_l\simeq (M_0^c\o
M_{\l_1}^c\o\dots\o M_{\l_n}^c)^{sing}[\l-2l]$. On the other hand, by
Lemma~2.5 this is isomorphic to $M_{\l_1}^c\o\dots \o
M_{\l_n}^c[\l-2l]$. Combining these isomorphisms, we get an
isomorphism $H_l\simeq M_{\l_1}^c\o\dots \o M_{\l_n}^c[\l-2l]$. It is
easy to see that this isomorphism is given by formula (2.6), and
satisfies the required properties.
\qed\enddemo

We will also need a modification of this theorem. As before, let
$z_0,\dots, z_n\in\R, z_0<\dots< z_n$, and let $\l_0=0$. Then we have
an embedding of the unions of hyperplanes: $\Cal C(z_1,\dots,
z_n)\subset \Cal C(z_0,z_1,\dots, z_n)$, which induces a map of homologies

$$i:H_l(z_1,\dots, z_n;\l_1,\dots,\l_n)\to H_l(z_0,\dots,
z_n;0,\l_1,\dots, \l_n).$$

\proclaim{Theorem 2.6} Denote for brevity $M^c=M_{\l_1}^c\o\dots \o
M_{\l_n}^c$. Then the following diagram is commutative: 

$$
\CD
(M^c)^{sing}[\l-2l]@>{\varphi}>>	
		H_l(z_1,\dots, z_n;\l_1,\dots, \l_n)\\ 
@VV{}V					@V{i}VV\\			
M^c[\l-2l]\simeq (M_0^c\o M^c)^{sing}[\l-2l]@>{\varphi}>> 
		H_l(z_0, \dots, z_n;0,\l_1,\dots, \l_n)
\endCD
\tag 2.9
$$
where the first vertical line is the natural embedding 
$(M^c)^{sing}[\l-2l]\subset M^c[\l-2l]$. 
\endproclaim

This theorem immediately implies the following corollaries, proof of
which is trivial:

\proclaim{Corollary 2.7} 
1. The map $i: H_l(z_1,\dots, z_n;\l_1,\dots,\l_n)\to H_l(z_0,\dots,
z_n;0,\l_1,\dots, \l_n)$ is an embedding.

2. Let $z_0=0$, and let $C\in H_l$ be a combinatorial cycle such that
$0$ is not an endpoint of any of the arcs. Then the corresponding
vector $\varphi^{-1}(C)\in (M_{\l_1}^c\o \dots\o M_{\l_n}^c)[\l-2l]$
is singular.
\endproclaim

\demo{Proof of Theorem 2.6} First,  note that the composition of maps 
$(M^c)^{sing}[\l-2l]\hookrightarrow M^c[\l-2l]\to (M_0^c\o
M^c)^{sing}[\l-2l]$, where the last arrow is the isomorphism of
Lemma~2.5, is given by $v\mapsto v_0^*\o v$, which immediately follows
from the definition. Similarly, the embedding
$M^c[\l-2l]\hookrightarrow (M^c\o M_0^c)^{sing}[\l-2l]$ is given by
$v\mapsto v\o v_0^*$.

Next, let $z_1<\dots<z_n<z_0$ (note this change of order!). 
Denote for brevity $\bold z=(z_1,\dots, z_n), \bl=(\l_1,\dots,
\l_n)$. Then we have the following commutative diagram:

$$\CD
(M^c)^{sing}[\l-2l]@>{\varphi}>>	
		H_l(\bold z;\bl)\\ 
@VV{}V					@V{}VV\\			
M^c[\l-2l]\simeq (M^c\o M_0^c)^{sing}[\l-2l]@>{\varphi}>> 
		H_l(\bold z, z_0;\bl, 0)
\endCD
\tag 2.10
$$

Indeed, for cycles $C\in H_l(\bold z,\bl)$ of the
form shown below, (2.10) immediately follows from the definition of
$\varphi$ in Theorem~2.4. On the other hand, these cycles span the
homology space $H_l(\bold z, \bl)$, which follows from Theorem~2.4. 

$$\epsfbox{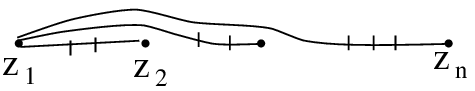}$$

Now we need to change the order of the points. Denote by $\check R$
the intertwining operator $PR_{M_0^c, M^c}:M_0^c\o M^c\to M_c \o
M_0^c$. It can be written as a product of $R$-matrices, each
interchanging $M_0^c$ with one of the factors $M^c_{\l_i}$. By
Theorem~2.4(2), we have the following commutative diagram:

$$\CD
(M^c\o M_0^c)^{sing}[\l-2l]@>{\varphi}>> 
		H_l(\bold z, z_0;\bl 0)\\
@VV{\check R^{-1}}V					@V{T^{-1}}VV\\
(M_0^c\o M^c)^{sing}[\l-2l]@>{\varphi}>> 
		H_l(z_0,\bold z;0, \bl)
\endCD
\tag 2.11
$$
 where the operator $T:H_l(z_0,\bold z;0,\bl)\to H_l(\bold z,
z_0;\l_1,\bl, 0)$ is the monodromy along the path shown on Figure~2.

\midinsert
\centerline{\epsfbox{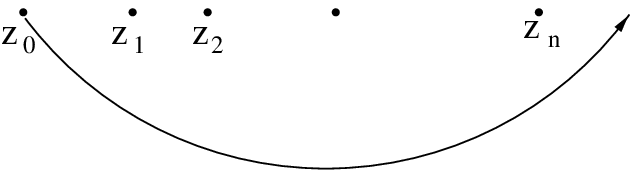}}
\botcaption{Figure 2}\endcaption
\endinsert

Combining (2.10) and (2.11), we get a commutative diagram as in the
statement of the theorem. It remains to check that composition of the
maps $(M^c)^{sing}\subset (M^c\o M_0^c)^{sing}@>{\check R^{-1}}>>
(M_0^c\o M^c)^{sing}$ coincides with the embedding
$(M^c)^{sing}\subset (M_o^c\o M^c)^{sing}$ in (2.9). This follows from
$\check R (v_0^*\o v)=v\o v_0^*$, which is
immediate from the definition of $\Cal R$. Similarly, the composition
of the maps

$$H_l(\bold z;\bl)\to H_l(\bold z,z_0;\bl, 0)@>{T^{-1}}>>
H_l(z_0,\bold z;0,\bl)$$ 
is equal to the map of homology spaces in (2.9): it suffices to check
this for a combinatorial cycle in $H_l(\bold z;\bl)$ in which case
it is obvious.
\qed\enddemo

\head 3. Complex conjugation and canonical basis \endhead

As before, we assume that we have fixed data $n,l,\kappa,
\bl=(\l_1,\dots, \l_n), \bold z=(z_1,\dots, z_n)$ as in the beginning
of Section~2. from now on, we assume in addition that $\l_i\in
\R$. Note also that since we assumed $\kappa\in \R$, we have $\bar
q=q^{-1}$. From now on, we will note need the dependence of the
homology space on the points $z_i$, so we will just write $H_l$ for
the homology space defined by (2.3). We will identify $H_l$ with the
weight subspace $M^c_{\l_1}\o\dots\o M^c_{\l_n}[\l-2l]$ using the
isomorphism of Theorem~2.3. 

Let $\overline{\phantom{T}}:\U\to \U$ be the antilinear algebra automorphism
defined by 

$$\overline{F}=F,\quad \overline{E}=E,\quad \overline{q^h}=q^{-h}
\tag 3.1$$

Define the associated antilinear involution on the  Verma  module $M_\l$
by $\overline{x v_\l}=\bar x v_\l, x\in \U$. More generally, if $M_1,
M_2$ are highest weight modules then define an  antilinear
involution $\psi:M_1\o M_2\to M_1\o M_2$ by 

$$\psi ( v\o v')=\overline{\Theta}  (\bar v\o \overline{v'}),
\tag 3.2$$
where $\Theta$ is defined by (1.6). One can check that $\psi^2=1$. 

\proclaim{Theorem 3.1 \rm (Lusztig)} Let $V, V'$ be irreducible
finite-dimensional $\U$-modules with the highest-weight vectors $v_\l,
v_\mu$ respectively. Let $b_i=F^{(i)} v_\l, b'_i=F^{(i)}
v_\mu$ be the bases in $V, V'$. Then there exists a unique
basis $b_{ij}=b_i\diamondsuit b'_j$  in $V\o V'$ such that:

\roster\item
	$b_i\diamondsuit b'_j=b_i\o b'_j+\sum_{k\ne 0} a_k b_{i-k}\o
	b'_{j+k}$, and $a_k\in q^{-1}\Z[q^{-1}]$
\item $\psi (b_i\diamondsuit b'_j)=b_i\diamondsuit b'_j$.
\endroster

Moreover, the coefficients $a_k$ in \rom{(1)} are non-zero only for $k>0$.
\endproclaim

The basis introduced above is called the canonical basis and has many
remarkable properties, which can be found in \cite{L}.

We will also need the dual basis. Let $M$ be a highest-weight module,
and $M^c$ the corresponding contragredient module. Define the
involution on $M^c$ by $\<\overline{v^*}, v\>=\overline{\<v^*, \bar
v\>}$. As before, let $b_k^*=(F^{(k)} v_\l)^*$ be the basis in $M^c$
dual to the basis $b_k=F^{(k)} v_\l$ in $M$. Since for $\l\in \Z_+$ we
have an embedding $V_\l\subset M_\l^c$, the elements $b_k^*, 0\le k\le
\l$ also give a basis in $V_\l$, dual to the basis $b_k$ with respect
to the contragredient form.

\proclaim{Theorem 3.2} Let $V, V'$ be irreducible finite-dimensional
$\U$-modules. Denote by $(b_i\diamondsuit b'_j)^*$ the basis in $V\o
V'$ dual to the canonical basis in $V\o V'$. Then this basis has
the following properties, which uniquely determine it:

\roster\item 

$$(b_i\diamondsuit b'_j)^*=b_i^* \o {b'_j}^*
	+\sum_{k>0}  \a_k  b_{i+k}^* \o {b'_{j-k}}^*,
	 \quad \a_k \in q^{-1}\Z[q^{-1}] \tag 3.3$$

\item 	$\psi^c((b_i\diamondsuit b'_j)^*)=(b_i\diamondsuit b'_j)^*$,
	where 
$$\psi^c(v^* \o (v')^*)= \tau(\Theta)
	( \overline{v^*}\o\overline{(v')^*}).\tag 3.4
$$
\endroster
Here $\Theta, C$ are defined by \rom{(1.6)}. 
\endproclaim
\demo{Proof} The theorem immediately follows from the definitions.
\qed \enddemo

From now on, we will call the basis, dual to the canonical basis,
``the dual canonical basis''. There are explicit formulas for the
canonical basis and the dual canonical basis -- see \cite{FK}.

More generally, one can define the canonical basis (and therefore, the
dual canonical basis) in a tensor product of any number of irreducible
finite-dimensional representations. This construction is also due to
Lusztig; we will use it in a more explicit form, taken from \cite{FK}. 

Let $V_1, \dots, V_n$ be irreducible finite-dimensional
representations. For every $i=1, \dots, n-1$ let $\check \Cal R_i:
V_1\o \dots\o V_n\to V_1\o\dots\o V_{i+1}\o V_i\o\dots\o V_n$ be the
operator $\check \Cal R=P\Cal R$ acting on $V_i\o V_{i+1}$. Let
$\sigma_0$ be the longest element in the symmetric group $\Sigma_n$:
$\sigma_0: (1 \dots n)\mapsto (n\dots 1)$, and let $\sigma_0=s_{i_1} \dots
s_{i_l}$ be its reduced expression. Define

$$\check \Cal R^{(n)}=\check \Cal R_{i_1}\dots \check \Cal R_{i_l}: 
	V_1\o\dots\o V_n \to V_n\o \dots \o V_1.
\tag 3.5$$

\proclaim{Proposition 3.3} 
\roster\item $\check \Cal R^{(n)}$ does not depend on the choice of
	the reduced expression for $\sigma_0$.
\item $\check \Cal R^{(n)}=\sigma_0 \Cal R^{(n)}$, where $\Cal R^{(n)}$
is defined inductively by $\Cal R^{(2)}=\Cal R$, 

$$\Cal R^{(n)}= (1\o \Cal R^{(n-1)}) \cdot (1\o \Delta^{n-2})(\Cal R)
	=(\Cal R^{(n-1)}\o 1)\cdot(\Delta^{n-2}\o 1)(\Cal R)
$$
and $\sigma_0: V_1\o\dots\o V_n \to V_n\o \dots \o V_1$ is the
permutation.  
\item Let 

$$\aligned
C^{(n)}=&(1\o C^{(n-1)})\cdot(1\o \Delta^{n-2})(C)
	=(C^{(n-1)}\o 1)\cdot(\Delta^{n-2}\o 1)(C)
	=q^{\half \sum_{i<j} h_i\o h_j},\\
\Theta^{(n)}=&(1\o \Theta^{(n-1)})\cdot(1\o \Delta^{n-2})(\Theta)
	=(\Theta^{(n-1)}\o 1)\cdot(\Delta^{n-2}\o 1)(\Theta)
\endaligned
$$
Then $\Cal R^{(n)}=C^{(n)}\Theta^{(n)}$.
\item 
$$\tau(\Theta^{(n)})
	=\check\Cal R ^{(n)}(C^{(n)})^{-1}\sigma_0.
$$
\endroster\endproclaim

\demo{Proof} (1) follows from the Yang-Baxter equation; (2) and (3)
are proved in \cite{FK, 1.2}. To prove (4), note that the  obvious
identities $\tau(\Cal R)=\Cal R^{21}, \tau(C)=C$ and the fact that
$\tau$ is a coalgebra automorphism imply $\tau(\Cal R^{(n)})=
\sigma_0(\Cal R^{(n)})=\sigma_0\cdot\Cal R^{(n)}\cdot \sigma_0, \tau
(C^{(n)}) =C^{(n)}$, after which we can use  $\Cal
R^{(n)}=C^{(n)}\Theta^{(n)}$. 
\qed\enddemo

Now we  define the dual canonical basis in a tensor product of $n$
representations: 

\proclaim{Proposition 3.4 \rm(see \cite{FK})} Let $V_1, \dots, V_n$ be irreducible finite
dimensional representations, and let $(F^{(\bold m)})^*$ be the dual
monomial basis \rom{(2.5)} in $V_1\o \dots\o V_n$. Then there exists a
unique basis $b_{\bold m}$ in $V_1\o \dots\o V_n$ such that

\roster \item 
$$b_{\bold m}=(F^{(\bold m)})^* +\sum_{\bold k\ne \bold m}a_{\bold k}
(F^{(\bold k)})^*, \quad a_{\bold k}\in q^{-1}\Z[q^{-1}].$$

\item $\psi^c(b_{\bold m})=b_{\bold m}$, where 

$$\psi^c(v_1\o\dots \o v_n)=
\tau(\Theta^{(n)})(\overline{v_1}\o\dots\o \overline{v_n}).\tag 3.6$$
\endroster
This basis is called the dual canonical basis.
\endproclaim

\remark{Remark} 
Similar to the case $n=2$, Lusztig defines a canonical basis in a
tensor product of $n$ modules (see \cite{L, Chapter~27}). One can
easily prove that the basis defined in Proposition~3.4 is dual to the
canonical basis of Lusztig. However, it is not needed for our
purposes.
\endremark

The main goal of our paper is to describe the dual canonical basis in
terms of the homology space $H_l$ defined by (2.3).

\proclaim{Lemma-Definition 3.5} Define an antilinear map $C_k(\C^l,
\Cal C,;\Cal S)\to C_k(\C^l, \Cal C;\Cal S)$ as follows: if $\Delta$
is a singular simplex, $s$ a section of the local system $\Cal S$ over
$\Delta$, then let $(\Delta, s)\mapsto (\bar\Delta, \overline{s(\bar
z)})$, where a bar denotes the standard complex conjugation in $\C^l$.
Then this map of complexes induces an antilinear involution on $H_l$
which will be denoted by $\overline{\phantom{T}}$.
\endproclaim

The first main result of this paper is the following theorem:

\proclaim{Theorem 3.6} Under the assumptions of the previous Lemma, the
isomorphism $M_{\l_1}^c\o \dots \o M^c_{\l_n}[\l-2l]\simeq H_l$
constructed in Theorem~\rom{2.3} identifies the complex conjugation in
$H_l$  with the involution $\psi^c$ defined by
\rom{(3.6)}
\endproclaim
\demo{Proof} It suffices to check it on the basis elements $(F^{(\bold
m)})^*$, where it is proved by direct calculation using Theorem~2.3 and
the identity $\tau(\Theta^{(n)})=\check\Cal R
^{(n)}(C^{(n)})^{-1}\sigma_0$ (see Proposition~3.3(4)).
\qed\enddemo

\remark{Remark 3.7} In fact, an obvious analog  of Theorem~3.6
holds for any simple Lie algebra $\g$. 
\endremark

\head 4. Explicit construction of the dual canonical basis. \endhead

In this section we give a geometric construction of the canonical
basis in a tensor product of irreducible
representations. This construction is a geometric counterpart of the
algebraic construction in \cite{FK}.

Let  $\l_1, \dots, \l_n\in \Z_+$. Denote for brevity

$$M=M_{\l_1}\o \dots\o M_{\l_n}, \quad V=V_{\l_1}\o \dots\o V_{\l_n}. 
\tag 4.1$$
Then we have an embedding 

$$V\subset M^c.\tag 4.2$$

We give a geometric construction of the dual canonical basis in the
weight subspace $V[\l-2l]\subset M^c[\l-2l]\simeq H_l$. Let $\l_i\in
\Z_+$ be such that $\sum \l_i\ge l$.  Let $B_l$ be the set of all
combinatorial cycles \rom{(}see Lemma-definition~\rom{2.1)} in the
space $H_l$ such that:

$$\aligned 
 &\text{\rom{(1)} On each arc, there is exactly one point $x$ placed.}\\ 
 &\text{\rom{(2)} All the arcs are in the upper half plane.}\\
 &\text{\rom{(3)} For each $i$, $z_i$ serves as an endpoint for at most
 	$\l_i$ arcs.}\\ 
 &\text{\rom{(4)} If for some $i$,  $z_i$ serves as an endpoint for less than
		$\l_i$ arcs}\\
	&\quad\text{then there are no arcs passing over $z_i$.} 
 \endaligned\tag 4.3$$

An example of a cycle $b\in B_l$ is given on Figure~3; here $n=6, l=7$
and $\l_1=3,\l_2=2,\l_3\ge 4, \l_4=\l_5=1, \l_6\ge 1$.
\midinsert
\centerline{\epsfbox{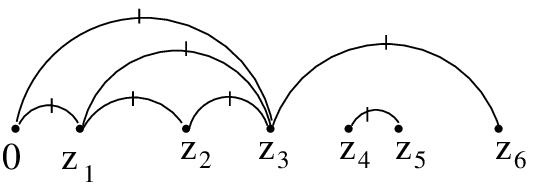}}
\botcaption{Figure 3} An example of a cycle $b\in B_l$.\endcaption 
\endinsert

\proclaim{Proposition 4.2} For every $b\in B_l$ define $\bold a(b)=(a_1,
\dots, a_n)\in \Z_+^n$ by the rule 
$a_i=(\text{number of arcs connecting
$z_i$ with all $z_j, j<i$})$. Then $b\mapsto \bold a(b)$ is a
bijection between $B_l$ and the set 

$$P_{\bl}(l)=\{\bold a\in \Z_+^n|a_i\le \l_i, \sum a_i=l\}
\tag 4.4$$
 
\endproclaim
\demo{Proof}
This proposition is of purely combinatorial nature, and can be proved
by induction on $l$ using the following easily verified fact: if we
remove from $b\in B_l$ the extreme right arc then we get a cycle from
$B_{l-1}$. Conversely, given a cycle $b\in B_l$ and an index $i, 1\le
i=le n$ there exists at most one way to add an arc to $b$ so that the
right end of this arc is $z_i$, and so that the resulting cycle lies
in $B_{l+1}$.  
\qed
\enddemo

This proposition allows us to index the elements of $B_l$
by $\bold m\in P_{\bl}(l)\subset \Z_+^n$: we will write
$b=b_{\bold m}$ if $b$ is an element of $B_l$ such that $\bold
a(b)=\bold m$.

\proclaim{Theorem 4.3} The elements $b_{\bold m}\in B_l, \bold m\in
P_{\bl}(l)$  form the dual canonical basis in $V[\l-2l]\subset
M^c[\l-2l]\simeq H_l$.

\endproclaim

Note that in the combinatorial construction of the dual canonical
basis in \cite{FK}, the basis is also parameterized by the pictures of
the form (4.3). However, the interpretation of these pictures is
different there.

The remaining part of this section and the two subsequent sections are
devoted to the proof of this theorem. The proof consists of several
steps.

\proclaim{Proposition 4.4} Every $b\in B_l$ 
lies in $V$. \endproclaim
\demo{Proof}
Consider the following finite collection of
infinite $l$-dimensional chains $C$ in $Y$: 

$$C_{ k_1, \dots, k_n}=\vcenter{\epsfbox{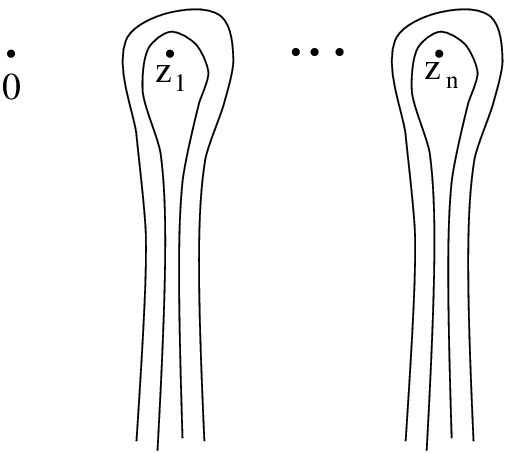}}
$$
where there are $k_1$ loops around $z_1$, $k_2$ loops around $z_2$, and
so on, with $\sum k_i=l$. Fix a section of the local system $\Cal S^*$
over each chain.  Then each such a chain $C$ defines a linear functional
on $H_l$ by the rule $f\mapsto \<f, C\>$, where $\<, \>$ is the intersection
pairing. One easily sees that this functional is well-defined.

 Due to Theorem~2.3, one can choose some sections over $C_{k_1,
\dots, k_n}$ so that the chains  $C_{k_1,\dots, k_n}$ will form a basis
dual to the monomial basis (2.5) in $M^c$. Thus, $v^*\in V\subset M^c $ if and
only if  $\<v^*, C_{ k_1,, \dots, k_n}\>=0$ for every $(k_1, \dots, k_n)$ such
that $k_i>\l_i$ for some $i$.

It follows from (4.3) that every $b\in B_l$  satisfies 
this condition and  thus, $b\in V$.
\qed\enddemo

The crucial step in proving Theorem~4.3 is proving the next two
propositions.

\proclaim{Proposition 4.5} 
For every $b=b_{\bold m}\in B$, 

$$b_{\bold m}=(F^{(\bold m)})^* 
	+\sum_{\bold k>\bold m} c_{\bold k}(F^{(\bold k)})^*,$$
where  $(F^{(\bold m)})^*$ is the monomial basis \rom{(2.5)}, 
$c_{\bold k}\in q^{-1}\Z [q^{-1}]$ and $<$ is the lexicographic
order on $\Z_+^n$. 
\endproclaim

\remark{Remark} In fact, it follows from  $b\in V$ that the sum above
can only contain monomials $(F^{(\bold k)})^*$ for $\bold k$ such that
$k_i\le \l_i, \sum k_i=l$. 
\endremark

\proclaim{Proposition 4.6} Every $b\in B_l$ is real,
i.e. $\overline{b}=b$. 
\endproclaim

The proofs of Propositions~4.5, 4.6  will be given in  Sections~5 and 6.

Now we can easily prove Theorem~4.3. Propositions~4.4, ~4.5 imply that
elements of $B$ lie in $V$ and are linearly independent. On the other
hand, it is obvious from Proposition~4.2 that the number of elements
in $B_l$ is equal to the dimension of $V[\l-2l]$.  Therefore, $B_l$ is
a basis in $V[\l-2l]$. Propositions~4.5 and~4.6 show that this basis
satisfies the definition of the dual canonical basis.
\qed

This theorem has many  corollaries.

\proclaim{Corollary 4.7} Assume that $\l_i\ge l$ for all $i$. Then the dual
canonical basis is given by

$$b_{\bold m}=[m_1]!\dots [m_n]!
	\vcenter{\epsfbox{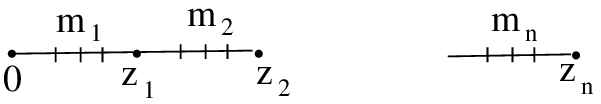}}
\tag 4.5$$
\endproclaim
\demo{Proof} Immediately follows from Theorem~4.3 and the following
lemma, which is easily proved by induction:

\proclaim{Lemma 4.8}
$$\vcenter{\epsfbox{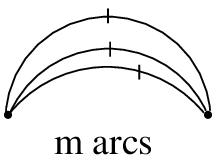}}=[m]!\vcenter{\epsfbox{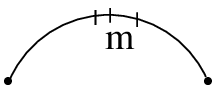}}
$$
\endproclaim
\qed\enddemo

Note that the cycles on the right-hand side of (4.5) are exactly the
anti-symmetric combinations of the bounded connected components of
$\R^l\setminus \Cal C_\R$. It has been known for a long time that
the connected components form a basis in the relative homology space
$H_L(\C^l, \Cal C; \Cal S)$ for any  local system $\Cal
S$ (see (\cite{A}, \cite{S} and references therein). Therefore, the
theorem above shows that the anti-symmetrization of this basis for the
local system described above  has a natural interpretation in
terms of  representation theory of $\U$: up to the factor $[m_1]!\dots
[m_n]!$, it is nothing else but the dual canonical basis in
$V_{\l_1}\o\dots\o V_{\l_n}$ under the assumptions that $\l_i\geq l$
for all $i$.

Another important special case is the dual canonical basis in the
subspace of singular vectors. Let us recall that a vector $v\in V$ is
called singular ($v\in V^{sing}$) if $ev=0$. It follows from the
general result of Lusztig (\cite{L, 27.2.5}) that there exists some subset
$B^{sing}\subset B$ of the dual canonical basis which is a basis in the
subspace $V^{sing}$. The following theorem describes this basis
explicitly.

\proclaim{Theorem 4.9} Let $B^{sing}$
be the subset of all cycles $b\in B_l$ satisfying the following
additional condition:

$$\text{There are no arcs ending at $0$.}
\tag 4.6$$

Then $B^{sing}$ is a basis in the subspace $V^{sing}\subset V$.
\endproclaim
\demo{Proof} By Corollary~2.7, we see that every $b\in B^{sing}$ is a
singular vector. Using Lemma~2.5 and simple combinatorial arguments,
similar to those used in the proof of Proposition~4.2, it is easy to
show that the dimension of $V^{sing}$ is equal to the number of
elements in $B^{sing}$.
\qed\enddemo

In particular, if $\l_1+\dots+\l_n-2l=0$ then Theorem~4.9 describes a
basis in the space $(V[0])^{sing}$, which coincides with the space of
invariants in $V$. Note that in this case condition (4.6) is
equivalent to the following condition:

$$\text{Each point $z_i$ is  an endpoint for exactly $\l_i$
arcs.}
\tag 4.7$$

Note that if we complete the real line by adding the point $\infty$ so
that the real line becomes a circle and the upper half-plane becomes a
disk then the cycles $b\in B^{inv}$ correspond to all the ways to join
$n$ points on a circle by non-intersecting arcs with multiplicities so
that the total number of arcs ending at $z_i$ is equal to
$\l_i$. Again, this parameterization of the basis in the subspace of
invariants is completely parallel to the combinatorial construction in
\cite{FK}. Both the geometric construction presented here and the
combinatorial construction in \cite{FK} make the cyclic symmetry of
the dual canonical basis in the space of invariants (see \cite{L,
28.2.9}) obvious.

\head 5. Proofs: the  case of 2-dimensional modules. \endhead

In this section we give the proofs of Propositions~4.5, 4.6 in the
special case where all $\l_i=1$, so that all $V_{\l_i}$ are
two-dimensional. Later we will show how to reduce the general case to
this special one.

Throughout this section, we assume $\l_i=1$ for all $i$. In this case,
conditions (4.3) imply that for every cycle $b\in B_l$ and every $i$,
the point $z_i$ is an endpoint of at most one arc. Here is a typical
example of a cycle $b \in B_l$:

$$\epsfbox{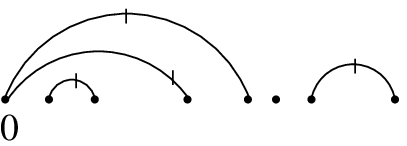}
\tag 5.1 $$

\proclaim{Lemma 5.1} In the notations of Section~\rom{2}, let $n=2$. Then
the cycle shown below represents the zero class in the anti-symmetrized
homology space $H_l$.

$$C=\vcenter{\epsfbox{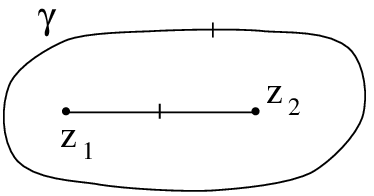}}
 $$
\endproclaim
\demo{Proof} Note first that it is indeed a cycle, since the monodromy
around $\gamma$ is equal to 1. Deforming this cycle so that
$\gamma$ approaches the interval $[z_1, z_2]$, we get that 

$$C=const(1+q^{-2}- 1 - q^{-2})\vcenter{\epsfbox{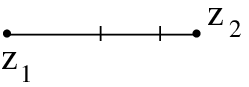}}
	=0$$
\qed\enddemo

\proclaim{Corollary 5.2} Let $A\in H_l$ be any combinatorial cycle satisfying
conditions \rom{(4.3, 4.6 )}. Then 

$$\epsfbox{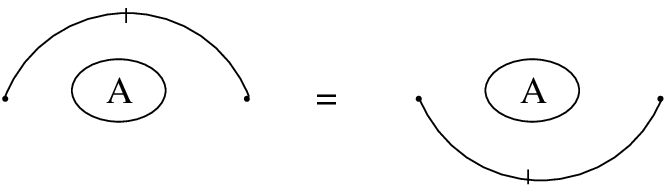}
\tag 5.2 $$
\endproclaim

Repeated application of this corollary yields the following result:

\proclaim{Proposition~5.3} Every cycle $b\in B_l$ 
satisfies $\bar b=b$. 
\endproclaim

This is nothing but Proposition~4.6 in our special case. 

\proclaim{Lemma 5.4} Let $A$ be any combinatorial cycle satisfying
conditions \rom{(4.3, 4.6)}, and
let $B, C$ be any combinatorial  cycles contained in the corresponding
regions \rom{(}dotted lines on the figure\rom{)}. Then 

$$\aligned
\vcenter{\epsfbox{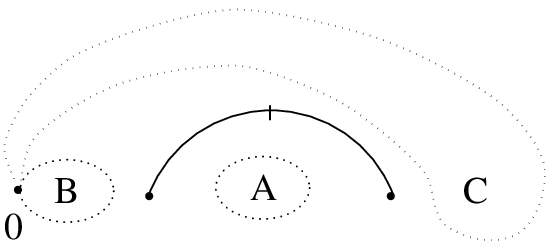}}=&\vcenter{\epsfbox{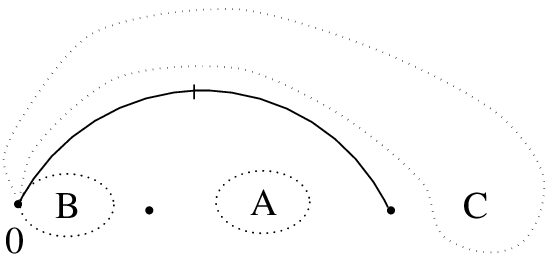}} \\
&-q^{-1}\vcenter{\epsfbox{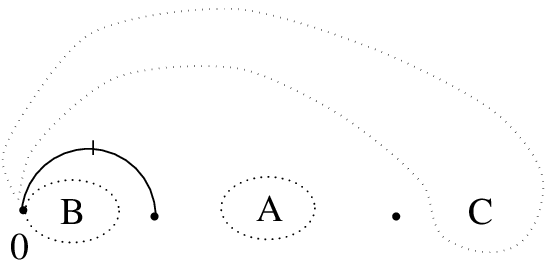}}
\endaligned
\tag 5.3 $$
\endproclaim

\demo{Proof} Obvious, since the ``total weight'' of the cycle $A$ is
equal to zero (that is, monodromy around the cycle encompassing $A$ is
trivial). \qed\enddemo

Using this proposition, we can prove by induction the following
corollary. 

\proclaim{Corollary 5.5} Let $b\in B_l$. Then we have the following
identity in $H_l\simeq M^c[\l-2l]$:

$$b=(F^{(\bold a(b))})^* +\sum_{\bold k>\bold a(b)} c_{\bold
	k}(F^{(\bold k)})^*,$$ where $c_{\bold k}\in q^{-1}\Z
	[q^{-1}]$, $<$ is the lexicographic order on $\Z^n$ and $\bold
	a(b)\in \Z_+^n$ is defined by $a_i=0$ if $z_i$ is the left
	endpoint of some arc, and $a_i=1$ otherwise.
\endproclaim

This is exactly  Proposition~4.5 in our special case, which
concludes the proof of Theorem~4.3 for the case $\l_i=1$.

\head 6. Proof: the general case
\endhead

In this section we will prove Propositions~4.5, 4.6 (and thus,
Theorem~4.3) for arbitrary $\l_i\in \Z_+$. We reduce the general
situation to the case $\l_i=1$, which has been proved already. Our
reduction is parallel to the approach in \cite{FK}.

Recall the setup of Section~3 and the anti-symmetrized homology space
$H_l$. 

Let us choose points $0<\zt_1<\dots<\zt_\l, \l=\sum\l_i$, and consider
the configuration of hyperplanes $\tilde \Cal C\subset \C^l$ defined
as in (2.1) with $z_i$ replaced by $\zt_p, 1\le p\le \l$ and $n$
replaced by $\l$. Let $\tilde \Cal S$ be the local system on
$\C^l\setminus\tilde\Cal C$ defined by (2.2) with the changes as above
and with $\tilde \l_p=1$. Denote by $\tilde H_l$ the corresponding
anti-symmetrized homology space.

\proclaim{Theorem 6.1} There exists a map $i_\bl:\tilde H_l\to H_l$ such
that 
\roster\item It preserves the complex conjugation \rom{(}see
Lemma-Definition~\rom{3.5)}

\item Let  $\pi:\{1, \dots, \l\}\to \{1, \dots, n\}$ be defined by 
$\pi(i)=1$ if $1\le i\le \l_1$, $\pi(i)=2$ if $\l_1+1\le i\le
\l_1+\l_2$, and so on. Let $\tilde C\in \tilde H_l$ be any
combinatorial cycle lying in the upper half-plane. Then: $i_\bl(\tilde
C)=0$ if $\tilde C$ contains a curve connecting $\zt_i$ with $\zt_j$ such
that $\pi(i)=\pi(j)$. Otherwise, $C=i_\bl(\tilde C)$ is also a
combinatorial cycle lying in the upper half-plane and defined as
follows: to each curve in $\tilde C$ connecting $\zt_i$ with $\zt_j$ there
corresponds a unique curve in $C$ connecting $z_{\pi(i)}$ with
$z_{\pi(j)}$, and with the same number of points $x$ marked on it. 

\endroster
\endproclaim

\example{Example} Let $n=2, \l_1=\l_2=2$. Then

$$\vcenter{\epsfbox{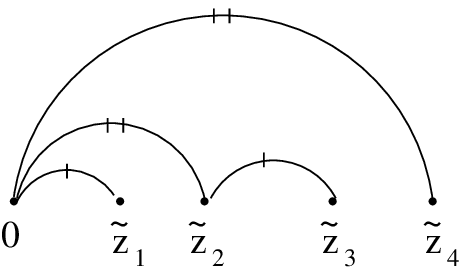}}@>{i_\bl}>>\vcenter{\epsfbox{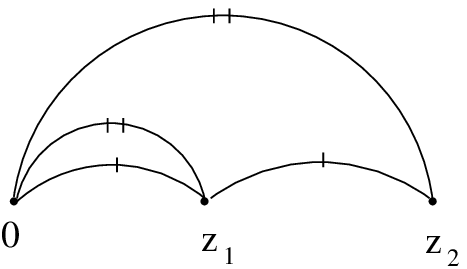}}
$$
\endexample
\demo{Proof}
Let us choose closed non-intersecting disks $D_1, \dots, D_n$ in $\C$
such that $D_i$ encloses the points $\zt_k$ with $\pi(k)=i$ (see
Figure~3). Choose a
continuous map $\phi:\C\to \C$ such that

\roster\item $\phi(D_i)=z_i, \phi(0)=0$

\item $\phi$ gives a diffeomorphism $\C\setminus\{D_1, \dots, D_n\}
	\simeq \C\setminus\{z_1, \dots, z_n\}$.

\item $\phi(\bar z)=\overline{\phi(z)}$
\endroster

\midinsert
\centerline{\epsfbox{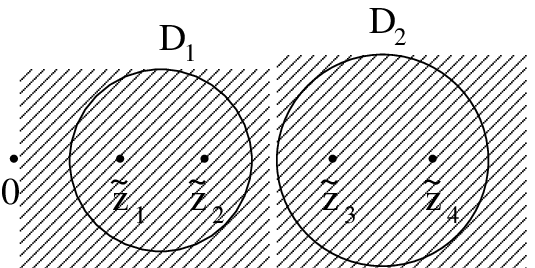}}
\botcaption{Figure 3.}\endcaption
\endinsert

It is easy to see that such a map exists. Let
$\phi^l=\phi\times\dots\times \phi:\C^l\to \C^l$. Then $\phi(\tilde
\Cal C)=\Cal C$, and thus $\phi$ induces a map of complexes of
relative chains $C_k(\C^l, \tilde\Cal C)\to C_k(\C^l, \Cal C)$. 

To extend it to chains with coefficients in local systems, we need one
more observation. Denote $U=\{\bold x\in \C^l|x_i \ne x_j, x_i\notin
D_k, x_i\ne 0\}$. The following lemma can be easily obtained by
comparing monodromies:

\proclaim{Lemma} The following local systems on $U$ are isomorphic:

$$\tilde\Cal S|_{U}\simeq ((\phi^l)^* \Cal S)|_U\tag 6.1$$

\endproclaim

Let us fix an isomorphism  in (6.1) by the condition that it identifies
$Br (\tilde \psi)$ with $Br(\psi)$ for $0\ll Re\ x_1\ll\dots\ll Re\ x_l$.

Now we can extend the map of relative chains to the relative chains with
coefficients in local systems: if $\tilde \Delta$ is a simplex in $\C^l$ and
$\tilde s$ is a section of $\tilde \Cal S$ over $\Delta\cap (\C^l\setminus
\tilde \Cal C)$ then let $\phi(\tilde \Delta, \tilde s)=(\Delta, s)$,
where $\Delta=\phi(\tilde\Delta)$ and $s$ is the image of $\tilde s$
under the morphism $\Gamma(\tilde\Delta\cap\C^l\setminus\tilde\Cal
C,\tilde\Cal S)\to  \Gamma(\tilde\Delta \cap U,\tilde\Cal S)\simeq 
\Gamma(\Delta\cap (\C^l\setminus \Cal C), \Cal S)$. One easily
verifies that this  map is indeed a morphism of complexes and
commutes with the action of the symmetric group $\Sigma_l$; thus, it gives
rise to a map of homologies $i_\bl:\tilde H_l\to H_l$. It is easy to
verify that so defined $i_\bl$ satisfies all the conditions of the
theorem. (This uses the special choice of sections of the local system
over combinatorial cycles, made in Definition~2.1.)
\qed\enddemo

\remark{Remark 6.2} It can be shown that if we rewrite the map $i_\bl$
using the identification $H_l\simeq M^c[\l-2l]$ \rom{(}see
Theorem~\rom{2.3)}, then it coincides with the dual to the canonical
embedding

$$ M_{\l_1}\o\dots\o M_{\l_n}\to 
M_1^{\o \l_1}\o M_1^{\o \l_2}\o \dots \o M_1^{\o \l_n}
$$

\endremark

\proclaim{Corollary 6.3} Propositions~\rom{4.5, 4.6} hold for arbitrary
$\l_i\in \Z_+$. 
\endproclaim
\demo{Proof} Using the map $i_\bl$ constructed in the previous theorem
we see that every $b\in B$ can be written as $i_\bl(\tilde b)$ for
some $\tilde b\in \tilde B_l$. Since $i_\bl$ commutes with the complex
conjugation, $\bar b=b$ follows from $\overline{\tilde b}=\tilde b$,
which was proved in the previous section. Similarly, the triangularity
condition follows from similar statement for $\tilde b$ (Corollary
5.5) and Lemma~4.8 allowing one to rewrite the monomial basis
$(F^{(\bold m)})^*$ in terms of the arcs with only one marked point.

\qed\enddemo

Therefore, we have proved Theorem~4.3, which gives us a construction
of the dual canonical basis in $V_{\l_1}\o\dots\o V_{\l_n}$ for arbitrary
highest weights $\l_i$.


\Refs
\widestnumber\key{AAA}

\ref\key A\by Aomoto, K.\paper Gauss-Manin connection of integral of
difference products
\jour J. Math. Soc. Japan \vol 39\yr 1987\pages 191--207
\endref

\ref\key FK\by Frenkel, I. and Khovanov, M.
\paper Canonical bases in tensor products and graphical calculus for
$U_q(\sltwo)$
\paperinfo preprint, August 1995
\endref

\ref\key KK\by Khovanov, M. and Kuperberg, G. 
\paper  Web bases for $sl(3)$ are not dual canonical
\paperinfo Preprint, available from {\sl
http://www.math.ucdavis.edu/\~{}greg/} 
\endref


\ref\key L\by Lusztig, G. \book Introduction to quantum groups\publ
Birkh\"auser \publaddr Boston \yr 1993\endref

\ref\key S\by Salvetti, M. \paper Topology of the complement of real
hyperplanes in $\C^n$ 
\jour Invent. Math.
\vol 88\yr 1987 \pages 603--618
\endref


\ref \key SV1\by Schechtman, V. and Varchenko, A.
\paper Arrangements of hyperplanes and Lie algebra homology
\jour Inv. Math.\vol 106 \yr 1991\pages 139--194
\endref



\ref\key SV2\bysame
\paper Quantum groups and homology of local systems
\inbook  Algebraic geometry and analytic geometry (Tokyo, 1990), 
 ICM-90 Satell. Conf. Proc.\pages 182--197
\publ Springer \publaddr Tokyo \yr 1991
\endref


\ref\key V1\by Varchenko, A. 
\book Multidimensional hypergeometric functions and representation
theory of quantum groups
\bookinfo Adv. Ser. Math. Phys. \vol 21
\publ World Scientific\publaddr River Edge, NJ
\yr 1995
\endref

\ref\key V2\bysame
\paper Asymptotic solutions to the Knizhnik--Zamolodchikov equation and
crystal base
\jour Comm. Math. Phys.\vol 171\yr 1995\pages 99--138\endref



\endRefs
\enddocument
\end